\documentclass{article}
\usepackage{authblk}
\usepackage{graphicx} 
\usepackage{amsmath}
\usepackage{amssymb}
\usepackage{amsthm}
\usepackage{amsmath}
\usepackage{subcaption}
\usepackage{float}
\usepackage{placeins}
\usepackage{mathrsfs}
\usepackage{dsfont}
\usepackage{bm}
\usepackage[table,xcdraw]{xcolor} 
\usepackage[numbib]{tocbibind}
\usepackage{enumerate}
\usepackage[normalem]{ulem}

\newtheorem{theorem}{Theorem}

\newtheorem{proposition}{Proposition}
\newtheorem{remark}{Remark}

\usepackage{stackengine,xcolor}
\usepackage{amsmath, amssymb}
\usepackage{mathrsfs}
\input pdf-trans

\allowdisplaybreaks

\newbox\qbox
\def\usecolor#1{\csname\string\color@#1\endcsname\space}
\newcommand\bordercolor[1]{\colsplit{1}{#1}}
\newcommand\fillcolor[1]{\colsplit{0}{#1}}
\newcommand\colsplit[2]{\colorlet{tmpcolor}{#2}\edef\tmp{\usecolor{tmpcolor}}%
	\def\tmpB{}\expandafter\colsplithelp\tmp\relax%
	\ifnum0=#1\relax\edef\fillcol{\tmpB}\else\edef\bordercol{\tmpC}\fi}
\def\colsplithelp#1#2 #3\relax{%
	\edef\tmpB{\tmpB#1#2 }%
	\ifnum `#1>`9\relax\def\tmpC{#3}\else\colsplithelp#3\relax\fi
}
\newcommand\outline[1]{\leavevmode%
	\def\maltext{#1}%
	\setbox\qbox=\hbox{\maltext}%
	\boxgs{Q q 2 Tr \thickness\space w \fillcol\space \bordercol\space}{}%
	\copy\qbox%
}

\bordercolor{white}
\fillcolor{black}
\def\thickness{.1}

\begin{document}


\title{An efficient scheme for approximating long-time dynamics of a class of non-linear models
}
		\author[a]{Jack Coleman
		}		
		\author[b]{Daozhi Han }
					
		\author [c,a]{Xiaoming Wang}
	
	        \affil[a]{Department of Mathematics and Statistics, Missouri University of Science and Technology, Rolla, MO 65409\\ jlcc7w@mst.edu}
		\affil[b]{Department of Mathematics, State University of New York at Buffalo, Buffalo, NY 14260\\ daozhiha@ub.edu}
		\affil[c]{Corresponding author, School of Mathematical Science, Eastern Institute of Technology, Ningbo, Zhejiang 315200, China\\ wxm@eitech.edu.cn}
		
		\date{\today}

\maketitle
%

\begin{abstract}
We propose a novel, highly efficient, mean-reverting-SAV-BDF2-based, long-time unconditionally stable numerical scheme for a class of finite-dimensional nonlinear models important in geophysical fluid dynamics. The scheme is highly efficient in that only a fixed symmetric positive definite linear problem (with varying right-hand sides) is solved at each time step. The solutions remain uniformly bounded for all time. We show that the scheme accurately captures the long-time dynamics of the underlying geophysical model, with the global attractors and invariant measures of the scheme converging to those of the original model as the step size approaches zero.

In our numerical experiments, we adopt an indirect approach, using statistics from long-time simulations to approximate the invariant measures. Our results suggest that the convergence rate of the long-term statistics, as a function of terminal time, is approximately first-order under the Jensen-Shannon metric and half-order under the total variation metric. This implies that extremely long simulations are required to achieve high-precision approximations of the invariant measure (or climate). Nevertheless, the second-order scheme significantly outperforms its first-order counterpart, requiring far less time to reach a small neighborhood of statistical equilibrium for a given step size.
\end{abstract}



\section{Introduction}

Many finite dimensional geophysical models are forced and damped together with a linear and a nonlinear term that conserves energy \cite{Majda2016}.  Well-known examples include the Lorenz 1996 model, Galerkin truncations of the damped-driven quasi-geostrophic models, etc. These models can be formulated using the following abstract format on a finite-dimensional Hilbert space $H$. 


\begin{equation}
	\frac{d\mathbf{u}}{dt} + A \mathbf{u} + \mathbf{N}(\mathbf{u}) = \mathbf{F}, \quad \mathbf{u} \in H, \quad  \mathbf{u} \mid_{t=0} = \mathbf{u}_0, \label{eq:Model}
\end{equation}
where the operators $A, N$ and the forcing term $\mathbf{F}$ satisfy the following assumptions
\begin{enumerate}[I.]
	\item $A\in B(H,H), A>0, \quad \exists \ell_0 > 0 \ s.t. \ ||A^{\frac{1}{2}} \mathbf{u}||^2 \geq \ell_0 ||\mathbf{u}||^2, \quad \forall \mathbf{u} \in H$.
	\item $\mathbf{N}(\mathbf{u})$ is {
		energy conservative} in the sense that $\mathbf{N}(\mathbf{u}) \cdot \mathbf{u} = 0, \quad \forall \mathbf{u} \in H$.
	\item $\mathbf{N}(\mathbf{u})$ is locally Lipschitz.
	\item $\mathbf{F} \in L^{\infty}(0, T, H), \quad ||\mathbf{F}||_{\infty} = \sup\limits_{t>0} ||\mathbf{F}(t)|| <\infty$.
\end{enumerate} 
Here $A$ represents the linear (symmetric) damping mechanism, $N$ represents the skew symmetric term (energy conservation term) that may contain linear and nonlinear effects.

The long time dynamics of this kind of models could be extremely complex with abundant chaotic or turbulent behavior, and numerical methods are indispensible in order to gain quantitative information.  For this type of systems, statistical properties are physically more relevant and robust \cite{SS2001, Majda2016, LM1994, MW2006, FMRT2001, MAG2005}.  In the case when the system is autonomous, the invariant measures, if exist, characterize the long time statistics, i.e., the climate, of the underlying model. Therefore, it is of great interest to develop numerical algorithms that are able to capture the long time statistics such as the invariant measure of these systems.  It was discovered that preserving the dissipativity of the underlying model is a key ingredient in designing numerical algorithms capable of capturing long time statistics through the so-called indirect approach  \cite{Wang2010, Wang2012, Wang2016}. See also \cite{E2001} and the references therein for the context of stochastic models.
Hence, the development of long-time stable and efficient method is a key component of designing algorithms that are capable of capturing long time statistics, especially since the indirect cost involves long-time simulation of the underlying system. 

There are at least two strategies for enhancing computational efficiency. The first is toemploy a higher-order in time method, allowing for relatively larger time steps can be taken with the same error tolerance. The second is to treat the nonlinear and skew symmetric terms explicitly, so that only linear, positive definite systems need to be solved at each time step. However, the latter approach typically leads to severe time step restrictions, as the explicit treatment can introduce instability. 
In a recent development, the authors of \cite{LSL2021} introduced a highly efficient BDF2-SAV-based IMEX scheme for solving the two-dimensional Navier-Stokes equations, which admits a uniform-in-time bound in a modified energy in the absence of external forcing. See also \cite{Yang2021b} for a closely related approach, termed ZEC. Nonetheless, the original BDF2-SAV scheme is non-autonomous, which is not conducive to dynamical system approach that we adopt for studying long time behavior. In addition, no boundedness of the solution is available for the case with non-zero external forcing, since small errors could accumulate in the auxiliary variable over long time. 
In order to suppress the growth of the error in the auxiliary variable over extended time, we introduce a simple mean-reverting mechanism in the scalar auxiliary variable.
Our novel {\bf mean-reverting-SAV-BDF2}-based IMEX scheme , where the auxiliary variable is governed by a simple mean-reverting equation is long-time stable. The scheme is extremely efficient, as it only involves solving a linear symmetric problem of the form of $(A+\omega I)\mathbf{u}=\mathbf{f}$ at each time step with the same $\omega$, but different $\mathbf{f}$. Our main result is that this novel mean-reverting-SAV-BDF2-based scheme is able to capture the long time statistics of the underlying model \eqref{eq:Model}. To the best of our knowledge, this is the first unconditionally stable scheme that treats the skew symmetric terms explicitly while maintaining uniform-in-time bound, without any time-step restriction in the presence of external forcing. 

We also investigate the performance of the scheme via simulation on long time intervals $[0, T], T\gg 1$, in order to approximate the invariant measure.
Our numerics suggests a half-order convergence rate to the invariant measure as a function of the terminal time $T$ using the total variation metric, while the rate of convergence in Jensen-Shannon metric is of order one. This implies that very long time simulations are needed in order to obtain two or three significant digits accuracy, highlighting  a key challenge associated with the study of climate and climate change of deterministic (finite dimensional) (perfect) climate models using indirect approach.
We also compared the performance of the second order scheme to that of its first order counterpart. Although the rate of convergence are roughly the same, the second order scheme requires significantly less time to enter a small neighborhood of the statistical equilibrium. This is a strong indicator of the superiority of the second order scheme over the first order one.

The rest of the paper is organized as follows. We recall a basic result on the original model as well as the auxiliary model in section 2.  The novel mean-reverting-SAV-BDF2 scheme is introduced in section 3. The long-time stability of the scheme is established in section 4.  The so-called asymptotic consistency  is established in section 5. The distance between the solutions of the scheme and the model is investigated in section 6. The convergence of the global attractors is then established in section  7. The convergence of the invariant measures is established in section 8. We apply the novel mean-reverting-SAV-BDF2 scheme to the five mode Lorenz 96 model and study its long time statistics in section 9. We offer concluding remarks at the end.

\section{Long time dynamics of the original and auxiliary model}

\subsection{Long time dynamics of the original model}
It is easy to show that the basic model enjoys a global attractor if the forcing term is time-independent. More specifically, we have
\begin{proposition}
	Model \eqref{eq:Model} processes a global attractor $\mathscr{ A}$ if $\mathbf{F}(t) \equiv \mathbf{F} \in H$. Moreover,  the attractor is contained in a ball with radius $\frac{||\mathbf{F}||_{\infty}}{\ell_0}$.
\end{proposition}

\begin{proof}
	The local existence of solutions to \eqref{eq:Model} follows from the linearity of A and the local Lipschitz property of the non-linear term $\mathbf{N}$. The global existence of the solution follows from taking the inner product of equation \eqref{eq:Model} with $\mathbf{u}$.
	\begin{align*}
		\frac{1}{2} \frac{d}{dt} ||\mathbf{u}(t)||^2 + ||A^{\frac{1}{2}} \mathbf{u}(t)||^2 = \mathbf{F}(t) \cdot \mathbf{u}(t) 
		\leq \frac{\ell_0}{2} ||\mathbf{u}(t)||^2 + \frac{1}{2\ell_0} ||\mathbf{F}||^2_{\infty}.
	\end{align*}
	Which implies 
	$$
	\frac{d}{dt} ||\mathbf{u}(t)||^2 + \ell_0||\mathbf{u}(t)||^2 \leq \frac{1}{\ell_0} ||\mathbf{F}||^2_{\infty} $$
	which further implies
	\begin{equation}     ||\mathbf{u}(t)||^2 \leq e^{-\ell_0 t} ||\mathbf{u}_0||^2 + \frac{1}{\ell_0^2} ||\mathbf{F}||^2_{\infty} \left( 1 - e^{-\ell_0 t} \right). \label{eq:Proof_Novel_SAV_u_bounded}
	\end{equation}
	The uniqueness follows from the local Lipschitz condition on $\mathbf{N}$. Notice that \eqref{eq:Proof_Novel_SAV_u_bounded} also implies that \eqref{eq:Model} possesses an attractive ball in $H$ with radius $(1+\delta)\frac{||\mathbf{F}||{\infty}}{\sqrt{\ell_0}}, \forall \delta>0$. Since any finite dimensional ball is pre-compact. We deduce, in the case when $\mathbf{F}$ is time independent, the existence of a global attractor following classical results. See for instance \cite{T1997}. 
\end{proof}

\subsection{The augmented system} 
We now introduce the following augmented model with an auxiliary scalar variable $q$.
\begin{subequations}
	\begin{align}
		&\frac{d\mathbf{u}}{dt} + A\mathbf{u} + q \mathbf{N}(\mathbf{u}) = \mathbf{F}, \label{eq:2_Novel_SAV_Model_u} \\
		&\frac{dq}{dt}  - \mathbf{N}(\mathbf{u}) \cdot \mathbf{u}= -\gamma q + \gamma,    \label{eq:2_Novel_SAV_Model_q} 
	\end{align}
	\label{eq:2_Novel_SAV_Model}
\end{subequations}
where the right-hand-side of \eqref{eq:2_Novel_SAV_Model_q} is mean-reverting, and $\gamma$ is a positive parameter at our disposal

It is easy to see that the last term on the left-hand-side of the scalar auxiliary variable identically vanishes due to the skew symmetry of $\mathbf{N}$. Hence, this augmented system reduces to the original model if $q(0)=1$.
Indeed, due to the mean-reverting nature of the auxiliary variable equation, $q$ approaches the value of 1 exponentially fast, leading to a system asymptotically equivalent to the original one at large time, regardless of the initial value of $q$.  This freedom of initial condition enables us to proceed with the dynamical system approach embedded in the Lax-type criteria \cite{Wang2010, Wang2016}.
It is also easy to see that the auxiliary system is dissipative and its global attractor is simply the product of the global attractor of the original model and the number 1.
\begin{proposition} 
	The augmented system \eqref{eq:2_Novel_SAV_Model} possesses a global attractor $\mathscr{A}_q = \mathscr{A} \times \{1\}$, when $ \mathbf{F}(t) \equiv \mathbf{F} \in H$.
\end{proposition}

\begin{proof} Straightforward energy method leads to 
	\begin{align}
		\frac{1}{2} \frac{d}{dt} \left( ||\mathbf{u}(t)||^2 + |q(t)|^2 \right) &+ ||A^{\frac{1}{2}} \mathbf{u}(t)||^2 + \gamma |q(t)|^2 \notag \\
		&= \mathbf{F}(t) \cdot \mathbf{u}(t) + \gamma q(t) \notag\\
		&\leq ||\mathbf{F}(t)|| \ ||\mathbf{u}(t)|| + \gamma |q(t)| \notag\\
		&\leq \frac{\ell_0}{2} ||\mathbf{u}(t)||^2 + \frac{1}{2\ell_0} ||\mathbf{F}(t)||^2  + \frac{\gamma}{2} |q(t)|^2 + \frac{\gamma}{2}.
	\end{align}
	
	Using the assumption that $\frac{1}{2}||A^{\frac{1}{2}} \mathbf{u}(t)||^2 \geq \frac{\ell_0}{2} ||\mathbf{u}(t)||^2$, we deduce 
	
	\begin{equation}
		\frac{1}{2} \frac{d}{dt} \left( ||\mathbf{u}(t)||^2 + |q(t)|^2 \right) + \frac{\ell_0}{2} ||\mathbf{u}(t)||^2 + \frac{\gamma}{2} |q(t)|^2 \leq \frac{1}{2\ell_0} ||\mathbf{F}(t)||^2 + \frac{\gamma}{2}.
	\end{equation}
	
	Denoting that $E(t) = ||\mathbf{u}(t)||^2 + |q(t)|^2$, we have
	
	\begin{equation}
		\frac{d}{dt} E(t) + \alpha E(t) \leq \frac{1}{\ell_0} ||\mathbf{F}||^2_{\infty} + \gamma ,
	\end{equation}
	
	where $$\alpha = \min\left\{ \ell_0, \gamma \right\} > 0, \quad ||\mathbf{F}||_{\infty} = \sup\limits_{t\geq 0} ||\mathbf{F}(t)||. $$
	
	Hence, 
	
	\begin{equation}
		E(t) \leq e^{-\alpha t} E(0) + \frac{1}{\alpha} \left( \frac{1}{\ell_0} ||\mathbf{F}||^2_{\infty} + \gamma  \right) \left( 1 - e^{-\alpha t} \right).
	\end{equation}
	
	This implies that \eqref{eq:2_Novel_SAV_Model} has a uniform-in-time bound and a global attractor $\mathscr{A}_{q}$ within a ball of radius $B_{R_0}$, with 
	\begin{equation}
		R_0^2 = \frac{1}{\alpha} \left( \frac{1}{\ell_0} ||\mathbf{F}||^2_{\infty} + \gamma \right). \label{eq:Gloabl_attract}
	\end{equation}
	It is also easy to see that  \eqref{eq:2_Novel_SAV_Model_q}  implies
	\begin{equation}
		|q(t) - 1|^2 \leq |q(0)-1|^2 e^{-\gamma t}.
	\end{equation}
	
	Hence $ \forall (\mathbf{u}, q) \in \mathscr{A}_q, \quad q = 1$. Denoting $\mathscr{A}$ the global attractor of  the original model. Then 
	\begin{equation}
		\mathscr{A}_q = \mathscr{A} \times \left\{1\right\}.
	\end{equation}
	
\end{proof}

\begin{remark} It is also easy to see that the set of invariant measures for the augmented system, denoted $\mathcal{IM}_q$, is intimately related to the set of invariant measures of the original model $\mathcal{IM}$. Indeed, one can check easily the following relationship.
$$\mathcal{IM}_q=\mathcal{M}\times \delta_q(q-1).$$
\end{remark}

\section{The mean-reverting-SAV-BDF2 scheme}

Here, we propose a second order efficient scheme that is based on BDF2 and SAV schemes. Our scheme is inspired by \cite{LSL2021}. {
	See also \cite{Yang2021} for a related approach under the name of ZEC.} However, there are at least two significant differences: (i) our auxiliary variable is mean-reverting  while no guarantee is provided for the numerical value of the auxiliary variable  to be close to the desired value of 1 at large time for the model in \cite{LSL2021}; and (ii) our scheme is autonomous while the model in \cite{LSL2021} is non-autonomous. 
{
	Our scheme also differs from ZEC or extended ZEC schemes in the sense that our mean-reverting mechanism is solution independent  vs solution dependent functions as presented in \cite{Yang2021, Yang2021b}. Our simple approach allows us to demonstrate that the auxiliary variable converges to the desired value of 1 as time approaches infinity, with an error bounded above by a constant multiple of the time-step.} We also point out that the SAV schemes proposed  in \cite{LSL2021} and  the one introduced in this paper differ greatly from the original SAV scheme designed for gradient flows \cite{Shen2018scalar} although they all involve a scalar auxiliary variable (SAV). We will use the mr-SAV in order to emphasize the mean-reverting nature of the auxiliary variable. The autonomous treatment is crucial for our dynamical system approach to the study of the long time behavior of the underlying model. It is also different from any other known scheme so far as we know. We shall demonstrate below that it can be used to approximate the long-time dynamics of the model in terms of approximating the attractor and invariant measures. 

Let $\gamma >0$ be a fixed parameter. Denote $\delta t = k >0$, the time-steps, we propose the following second order mr-SAV-BDF2 scheme. 

\begin{subequations}
	\begin{equation}
		\frac{3\mathbf{u}^{n+1} - 4\mathbf{u}^{n} + \mathbf{u}^{n-1}}{2 \delta t} + A\mathbf{u}^{n+1} + q^{n+1} \mathbf{N}(2\mathbf{u}^n - \mathbf{u}^{n-1}) = \mathbf{F}^{n+1}, \label{eq:2_Novel_SAV_discretization_u}
	\end{equation}
	\begin{equation}
		\frac{3q^{n+1} - 4q^{n} + q^{n-1}}{2 \delta t} + \gamma q^{n+1} - \mathbf{N}(2\mathbf{u}^n - \mathbf{u}^{n-1}) \cdot \mathbf{u}^{n+1} = \gamma , \label{eq:2_Novel_SAV_discretization_q}
	\end{equation}
	\label{eq:2_Novel_SAV_discretization}
\end{subequations}
where $\mathbf{u}^{n} \approx \mathbf{u}(n\delta t), \quad q^{n} \approx q(n\delta t) \approx 1, \quad q^0 = q^1 = 1$. 

Notice that \eqref{eq:2_Novel_SAV_discretization} is the discretization of \eqref{eq:2_Novel_SAV_Model}.

With the aid of $\left[ \frac{3}{2} I+ \delta t A \right]^{-1} $, the novel mr-SAV-BDF2 scheme can be solved via the following 
\begin{subequations}
	\begin{align}
		q^{n+1} &= \frac{1}{B^n} \left[ \frac{3}{2} + \delta t \gamma \right]^{-1} \left( \frac{4q^n - q^{n-1}}{2} + \delta t \gamma \right. \notag \\
		&+ \left. \delta t \mathbf{N}(2\mathbf{u}^{n} - \mathbf{u}^{n-1})  \cdot \left[ \frac{3}{2}I + \delta t A \right]^{-1} \left( \frac{4\mathbf{u}^{n} - \mathbf{u}^{n-1}}{2} + \delta t \mathbf{F}^{n+1} \right) \right),  \label{eq:simplified_Novel_SAV_q}\\
		\mathbf{u}^{n+1} &= \left[ \frac{3}{2} I + \delta t A \right]^{-1} \left(\frac{4\mathbf{u}^{n} - \mathbf{u}^{n-1}}{2} - \delta t q^{n+1} \mathbf{N}(2\mathbf{u}^{n} - \mathbf{u}^{n-1}) + \delta t \mathbf{F}^{n+1} \right), \\
		&\text{where } \notag \\
		B^n &= 1 + \delta t^2 \left[ \frac{3}{2} + \delta t \gamma \right]^{-1} \left( \left[ \frac{3}{2}I + \delta t A \right]^{-1} \mathbf{N}(2\mathbf{u}^{n} - \mathbf{u}^{n-1}) \right) \cdot \mathbf{N}(2\mathbf{u}^{n} - \mathbf{u}^{n-1}).
	\end{align}
\end{subequations}
This implies the super efficiency of the scheme as it only involves solving $( \frac{3}{2}I + \delta t A)U = F$ at each time step. One could even argue that the scheme is almost as efficient as it could be since the implicit treatment of the dissipative term is required for stability reason. 

\begin{remark} More general second order time discretization such as those A-stable generalized BDF2 schemes considered in \cite{WY2024} can be utilized in lieu of the classical BDF2. The long time stability analysis carries through by utilizing the analysis for the generalized BDF2 schemes presented in \cite{WY2024}. 
\end{remark}

\section{Long-time stability of the scheme}

Here we show that the solution to \eqref{eq:2_Novel_SAV_discretization} remains bounded for all time (all $n$) regardless of the initial conditions. We emphasize that $\mathbf{u}^1$, $q^1$ can take on arbitrary values, independent of the original initial data $\mathbf{u}^0, q^0$. This is different from classical numerical analysis. This added freedom in initial data is crucial to our dynamical system approach to approximating the long time statistics.

For this purpose we recall the classical G-matrix and G-norm associated with BDF2 scheme, see for instance \cite{HW2002},
\begin{equation}
	G = \frac{1}{4} \begin{bmatrix}
		1 & -2 \\
		-2 & 5 \\
	\end{bmatrix}
\end{equation}

Notice that G is positive-definite and symmetric. It induces a norm on $\mathbb{H} = H \times H$ and $\mathbb{R}^2$ via
\begin{subequations}
	\begin{align}
		||\mathbf{V}||^2_G &= \mathbf{V} \cdot G\mathbf{V}, \quad \mathbf{V} = \begin{bmatrix}
			\mathbf{v}_1 \\
			\mathbf{v}_2 
		\end{bmatrix} \in \mathbb{H} = H \times H, \\
		||\mathbf{Q}||^2_G &= \mathbf{Q} \cdot G\mathbf{Q}, \quad \mathbf{Q} = \begin{bmatrix}
			q_1 \\
			q_2 
		\end{bmatrix} \in \mathbb{R}^2,
	\end{align}
\end{subequations}
where the $\cdot$ denotes the inner product in corresponding Hilbert spaces.

Notice that the G-norm is an equivalent norm on $\mathbb{H}$ and $\mathbb{R}^2$. In fact, there exists $C_{l}, C_u >0$, s.t.

\begin{subequations}
	\begin{equation}
		C_{l} ||\mathbf{V}||^2 \leq ||\mathbf{V}||^2_G \leq C_u ||\mathbf{V}||^2, \quad \forall \mathbf{V} \in \mathbb{H},
	\end{equation}
	\begin{equation}
		C_{l} ||\mathbf{Q}||^2 \leq ||\mathbf{Q}||^2_G \leq C_u ||\mathbf{Q}||^2, \quad \forall \mathbf{Q} \in \mathbb{R}^2.
	\end{equation}
\end{subequations}

Denoting $\mathbf{V}^{n+1} = \begin{bmatrix}
	\mathbf{v}^{n} \\
	\mathbf{v}^{n+1}
\end{bmatrix}, \quad \mathbf{Q}^{n+1} = \begin{bmatrix}
	q^{n} \\
	q^{n+1}
\end{bmatrix}$.\\

$\eqref{eq:2_Novel_SAV_discretization_u} \cdot \mathbf{u}^{n+1} + \eqref{eq:2_Novel_SAV_discretization_q} \cdot q^{n+1}$  leads to, see for instance \cite{HW2002}

\begin{align*}
	&||\mathbf{V}^{n+1}||_G^2 + ||\mathbf{Q}^{n+1}||_G^2 - ||\mathbf{V}^{n}||_G^2 - ||\mathbf{Q}^{n}||_G^2 + \delta t ||A^{\frac{1}{2}} \mathbf{u}^{n+1}||^2 + \delta t |q^{n+1}|^2 \\
	&+ \frac{1}{4} || \mathbf{u}^{n+1} - 2\mathbf{u}^{n} +\mathbf{u}^{n-1} ||^2 + \frac{1}{4} | q^{n+1} - 2q^{n} + q^{n-1} |^2 \\
	&= \delta t \mathbf{F}^{n+1} + \gamma \delta t q^{n+1} \\
	&\leq \frac{\ell_0}{2} \delta t ||\mathbf{u}^{n+1}||^2 + \frac{\delta t}{2\ell_0} ||\mathbf{F}||^2_{\infty} + \frac{\gamma \delta t}{2} |q^{n+1}|^2 + \frac{\gamma \delta t}{2}.
\end{align*}

This  implies 

\begin{align}
	&||\mathbf{V}^{n+1}||_G^2 + ||\mathbf{Q}^{n+1}||_G^2 + \frac{\ell_0}{2} \delta t ||\mathbf{u}^{n+1}||^2  + \frac{\gamma}{2} \delta t |q^{n+1}|^2 \notag\\
	&+ \frac{1}{4} \left( ||\mathbf{u}^{n+1} - 2\mathbf{u}^{n} +\mathbf{u}^{n-1}||^2 + |q^{n+1} - 2q^{n} + q^{n-1}|^2\right) \notag\\
	&\leq ||\mathbf{V}^{n}||_G^2 + ||\mathbf{Q}^{n}||_G^2 + \frac{\delta t}{2\ell_0} ||\mathbf{F}||^2_{\infty} + \frac{\gamma \delta t}{2}. \label{eq:2_Novel_SAV_VQ_boundedness}
\end{align}

Let \begin{equation}
	\alpha = \frac{1}{6} \min \left\{ \ell_0, \gamma \right\} > 0,
\end{equation}
and define 
\begin{equation}
	E^n = ||\mathbf{V}^{n}||_G^2 + \alpha \delta t ||\mathbf{u}^n||^2 + ||\mathbf{Q}^{n}||_G^2 + \alpha \delta t |q^n|^2 .
\end{equation}

We deduce, after ignoring the last positive terms on the LHS of \eqref{eq:2_Novel_SAV_VQ_boundedness} and adding $\alpha \delta t \left( ||\mathbf{u}^{n}||^2 + |q^{n}|^2\right)$ to both sides of \eqref{eq:2_Novel_SAV_VQ_boundedness},
$$E^{n+1} + 3\alpha \delta t \left( ||\mathbf{u}^{n+1}||^2 + |q^{n+1}|^2\right) + \alpha \delta t \left( ||\mathbf{u}^{n}||^2 + |q^{n}|^2\right) \leq E^n + \frac{\delta t}{2\ell_0} ||\mathbf{F}||^2_{\infty} + \frac{\gamma \delta t}{2}.$$

Assume $\delta t \leq 1,$  and define $ \beta = \min \left\{ \alpha, \alpha C_{l} \right\}> 0,$ we have
\begin{subequations}
	\begin{equation}
		\left(1 + \beta \delta t\right)E^{n+1} \leq E^n + \left( \frac{1}{2\ell_0} ||\mathbf{F}||^2_{\infty} + \frac{\gamma}{2} \right) \delta t,
	\end{equation}
\end{subequations}

This implies 

\begin{equation}
	E^{n+1} \leq  \frac{1}{\left( 1 + \beta \delta t\right)^{n}}  E^{1} + \frac{1}{2\beta} \left( \frac{1}{\ell_0} ||\mathbf{F}||^2_{\infty} + \gamma \right) . \label{eq:2_Novel_SAV_uniform_bounded}
\end{equation}

Notice that the scheme \eqref{eq:2_Novel_SAV_discretization} generates  a discrete dynamical system on the product space $\mathbb{H} \times \mathbb{R}^2$ with the following solution semigroup
\begin{equation}
	\mathbb{S}_k \left(\begin{bmatrix}
		\mathbf{u}_0 \\
		\mathbf{u}_1
	\end{bmatrix}, \begin{bmatrix}
		q_0 \\
		q_1
	\end{bmatrix} \right) = \left(\begin{bmatrix}
		\mathbf{u}_1 \\
		\mathbf{u}_2
	\end{bmatrix}, \begin{bmatrix}
		q_1 \\
		q_2
	\end{bmatrix} \right).
\end{equation}
Hence, we have 

\begin{proposition}
	The second order scheme \eqref{eq:2_Novel_SAV_discretization} is uniformly bounded in the sense that \eqref{eq:2_Novel_SAV_uniform_bounded} holds for all initial data. Moreover,  \eqref{eq:2_Novel_SAV_discretization} generates a discrete disapatve dynamical system $\mathbb{S}_k$ on $\mathbb{H} \times \mathbb{R}^2$ with global attractor $\mathscr{A}_k\subset B_{R_1}(\mathbb{H})\times B_{R_1}(\mathbb{R}^2)$ where  $R_1 = \frac{1}{\sqrt{\beta \ell_0} }\sqrt{ ||\mathbf{F}||^2_{\infty} + \gamma\ell_0}$  if $\mathbf{F}$ is time dependent using the G-metric. 
\end{proposition}

\section{Asymptotic consistency of the system}

We have allowed arbitrary initial data on the product space in the scheme in order to conform to the dynamical system formalism. However, this also leads to the question of whether the solution to the scheme with arbitrary initial data converge to those of the underlying (augmented) model. 
However, if we focus on the global attractor $\mathscr[1]{A}_k$, 
we can show that the first and second component are always close,  with the difference bounded by a constant times $k$. 
This implies that points on the global attractors of the scheme are close to the initial conditions that we use for multi-step numerical methods. This further implies that the solution on the global attractor of $\mathbb{S}_k$ are close to solution to the original system after extension to the product space. This plays a key role in the convergence of the long time properties that we shall establish in the sequel.

{
	For notion purposes, the superscript $n$ denotes time-step count index, while the subscript $j, j=1,2$ denotes the first and second component of a a vector in the product space. For instance, $\mathbf{u}_j$ is the $j^{th}$ component of $\mathbf{V}=\begin{bmatrix}
		\mathbf{u}_1 \\
		\mathbf{u}_2
	\end{bmatrix}$.}

\begin{proposition}
	There exists a constant $C > 0$, independent of \(k\), s.t. 
	\begin{equation}
		||\mathbf{u}_1 - \mathbf{u}_2|| + |q_1 - q_2| \leq Ck, \quad \forall \left( \begin{bmatrix}
			\mathbf{u}_1 \\
			\mathbf{u}_2
		\end{bmatrix}, \begin{bmatrix}
			q_1 \\
			q_2
		\end{bmatrix} \right) \in \mathscr{A}_k, \forall k.
		\label{Eq:proposition_3}
	\end{equation}
\end{proposition}

\begin{proof}
	\eqref{eq:2_Novel_SAV_discretization} implies 
	\begin{subequations}
		\begin{equation*}
			||\mathbf{u}^{n+1} - \mathbf{u}^{n}|| \leq \frac{1}{3} ||\mathbf{u}^{n} - \mathbf{u}^{n-1}|| + \frac{2}{3}k\left( ||A\mathbf{u}^{n+1}|| + |q^{n+1}| \ ||\mathbf{N}(2\mathbf{u}^{n} -\mathbf{u}^{n-1})|| + ||\mathbf{F}||^2_{\infty} \right),
		\end{equation*}
		\begin{equation*}
			|q^{n+1} - q^n| \leq \frac{1}{3} |q^{n} - q^{n-1}| + \frac{2}{3}k\left( \gamma |q^{n+1}| + \gamma + |\mathbf{N}(2\mathbf{u}^{n} -\mathbf{u}^{n-1}) \cdot \mathbf{u}^{n+1}| \right).
		\end{equation*}
	\end{subequations}
	Recall from proposition 2 
	\begin{equation}
		||\mathbf{V}||^2 + ||\mathbf{Q}||^2 \leq \frac{C_u}{2\beta} \left( \frac{1}{\ell_0} ||\mathbf{F}||^2_{\infty} + \gamma \right) , \quad \forall (\mathbf{V}, \mathbf{Q}) \in \mathscr{A}_k, \quad \forall k.
	\end{equation}
	This, when combined with the local Lipschitz continuity of the nonlinear term $\mathbf{N}$, we have for a generic constant C, independent of k, s.t. 
	\begin{subequations}
		\begin{equation*}
			||\mathbf{u}^{n+1} - \mathbf{u}^{n}|| \leq \frac{1}{3} ||\mathbf{u}^{n} - \mathbf{u}^{n-1}|| + Ck,
		\end{equation*}
		\begin{equation*}
			|q^{n+1} - q^n| \leq \frac{1}{3} |q^{n} - q^{n-1}| + Ck,
		\end{equation*}
	\end{subequations}
	$\forall \left( \begin{bmatrix}
		\mathbf{u}_1 \\
		\mathbf{u}_2
	\end{bmatrix}, \begin{bmatrix}
		q_1 \\
		q_2
	\end{bmatrix} \right) \in \mathscr{A}_k, \quad \forall n \geq 1.$
	
	This implies that after iterating over n, 
	
	\begin{subequations}
		\begin{equation*}
			||\mathbf{u}^{n+1} - \mathbf{u}^{n}|| \leq \frac{1}{3^n} ||\mathbf{u}^{1} - \mathbf{u}^{0}|| + \frac{3}{2}Ck,
		\end{equation*}
		\begin{equation*}
			|q^{n+1} - q^n| \leq \frac{1}{3^n} |q^{1} - q^{0}| + \frac{3}{2} Ck.
		\end{equation*}
	\end{subequations}
	
	Since $\mathscr{A}_k$ is invariant, therefore $\forall \left( \begin{bmatrix}
		\mathbf{u}_1 \\
		\mathbf{u}_2
	\end{bmatrix}, \begin{bmatrix}
		q_1 \\
		q_2
	\end{bmatrix} \right) \in \mathscr{A}_k$ and $\forall n \in \mathbb{Z}^{+}, \exists \left( \begin{bmatrix}
		\mathbf{u}_n^0 \\
		\mathbf{u}_n^1
	\end{bmatrix}, \begin{bmatrix}
		q_n^0 \\
		q_n^1
	\end{bmatrix} \right) \in \mathscr{A}_k,$ s.t.  $ \mathbb{S}_k^n \left( \begin{bmatrix}
		\mathbf{u}_n^0 \\
		\mathbf{u}_n^1
	\end{bmatrix}, \begin{bmatrix}
		q_n^0 \\
		q_n^1
	\end{bmatrix} \right) = \left( \begin{bmatrix}
		\mathbf{u}_1 \\
		\mathbf{u}_2
	\end{bmatrix}, \begin{bmatrix}
		q_1 \\
		q_2
	\end{bmatrix} \right).$
	
	Hence, 
	\begin{subequations}
		\begin{equation*}
			||\mathbf{u}_1 - \mathbf{u}_2|| \leq \frac{1}{3^n} \|\mathbf{u}_n^1 - \mathbf{u}_n^0\| + \frac{3}{2} Ck,
		\end{equation*}
		\begin{equation*}
			|q_1 - q_2| \leq \frac{1}{3^n} |q^0_n - q^1_n| + \frac{3}{2} Ck, \quad \forall n.
		\end{equation*}
	\end{subequations}
	
	Letting $n$ approach $\infty$, we derive the desired result  \eqref{Eq:proposition_3} with a modified C which is independent of k.
\end{proof}

\begin{remark}
	An immediate consequence of the asymptotic consistency of the scheme is that the auxiliary variable approaches the desired value of 1 as the timestep index $n$ approaches infinity, with an error of order $\delta t = k$. To see this, we note that the nonlinear term in the scheme for the auxiliary variable \eqref{eq:2_Novel_SAV_discretization_q} satisfies
	$$\mathbf{N}(2\mathbf{u}^n - \mathbf{u}^{n-1}) \cdot \mathbf{u}^{n+1}
	=\mathbf{N}(2\mathbf{u}^n - \mathbf{u}^{n-1}) \cdot (\mathbf{u}^{n+1}-2\mathbf{u}^n + \mathbf{u}^{n-1}) = \mathcal{O}(k)$$
	for large $n$, thanks to the energy-conservative property of $\mathbf{N}$ and the asymptotic consistency.
	Hence,
	$$\frac{3q^{n+1} - 4q^{n} + q^{n-1}}{2 \delta t} + \gamma q^{n+1} -\gamma = \mathcal{O}(k)$$
	for large $n$, which implies
	$$\lim_{n\rightarrow\infty} |q^n-1| \le C k.$$
	A uniform-in-$n$ bound is also available, with dependence on the initial data included.
	This demonstrates that our scheme for the extended system is very close to the following efficient BDF2-Gear's extrapolation IMEX scheme for the original system:
	$$\frac{3\mathbf{u}^{n+1} - 4\mathbf{u}^{n} + \mathbf{u}^{n-1}}{2 \delta t} + A\mathbf{u}^{n+1} +  \mathbf{N}(2\mathbf{u}^n - \mathbf{u}^{n-1}) = \mathbf{F}^{n+1}.$$
	Although the solution to this scheme is not necessarily uniform-in-time bounded if the timestep is not small enough, our modified scheme mimics this scheme while achieving boundedness without any severe time step restriction. If the system is a spatial discretization of the incompressible Navier-Stokes equations, the modification of the Reynolds number due to the introduction of the auxiliary variable is small.
	Notably, such a bound on $q^n - 1$ is unavailable for any other schemes in the literature. Hence, it is not known a priori whether the alteration to the Reynolds number remains small for all time, even in this finite-dimensional case for those schemes. Therefore, our scheme is theoretically preferable for long-time turbulence studies.
\end{remark}

\section{Distance between discrete and continuous trajectories}

Here we consider the distance between the discrete trajectory $\mathbb{S}_k^n (\mathbf{V}^0, Q^0)$ and the continuous trajectory $\mathbb{S}(t) (\mathbf{V}^0, \mathbf{Q}^0)$, for $(\mathbf{V}^0, \mathbf{Q}^0) \in \mathcal{A}_k$ over finite time interval.  This is needed for the investigation of long time behavior.

Recall that 
$$\mathbb{S}(t) (\mathbf{V}^0, \mathbf{Q}^0) = \left( S_q(t) \begin{bmatrix}
	\mathbf{v}^0\\
	q^0
\end{bmatrix}, S_q(t)\begin{bmatrix}
	\mathbf{v}^1\\
	q^1 
\end{bmatrix} \right)
= \left( \begin{bmatrix}
	\mathbf{u}_0(t)\\
	q_0(t)
\end{bmatrix}, \begin{bmatrix}
	\mathbf{u}_1(t)\\
	q_ 1(t)
\end{bmatrix} \right)$$
here $S_q(t)$ is the solution semi-group associated with the auxiliary system \eqref{eq:2_Novel_SAV_discretization}.

Now assume that $\mathbf{F}$ and $\mathbf{N}$ are smooth enough so that the solutions to \eqref{eq:2_Novel_SAV_discretization} are sufficiently smooth. This implies, for $ (\mathbf{V}^0, \mathbf{Q}^0) \in B_{R_1}, \quad \exists C = C(k),$ s.t. 

\begin{subequations}
	\begin{align}
		&\frac{3 \mathbf{u}_0\left( (n+1)k\right) - 4 \mathbf{u}_0\left( nk\right) + \mathbf{u}_0\left( (n-1)k\right)}{2k} + A \mathbf{u}_0\left( (n+1)k\right) \notag \\
		&+ q_0\left( (n+1)k\right) \mathbf{N}(2\mathbf{u}_0(nk) -  \mathbf{u}_0\left( (n-1)k\right)) = \mathbf{F}((n+1)k) + \varepsilon_{u_0, n},
	\end{align}
	\begin{align}
		&\frac{3 q_0\left( (n+1)k\right) - 4 q_0\left( nk\right) + q_0\left( (n-1)k\right)}{2k} + \gamma q_0\left( (n+1)k\right) \notag \\
		&- \mathbf{N}(2\mathbf{u}_0(nk) -  \mathbf{u}_0\left( (n-1)k\right)) \cdot \mathbf{u}_0\left( (n+1)k\right) = \gamma + \varepsilon_{q_0, n},
	\end{align}
	
	\begin{align}
		&\frac{3 \mathbf{u}_1\left( (n+1)k\right) - 4 \mathbf{u}_1\left( nk\right) + \mathbf{u}_1\left( (n-1)k\right)}{2k} + A \mathbf{u}_1\left( (n+1)k\right) \notag \\
		&+ q_1\left( (n+1)k\right) \mathbf{N}(2\mathbf{u}_1(nk) -  \mathbf{u}_1\left( (n-1)k\right)) = \mathbf{F}((n+1)k) + \varepsilon_{u_1, n},
	\end{align}
	\begin{align}
		&\frac{3 q_1\left( (n+1)k\right) - 4 q_1\left( nk\right) + q_1\left( (n-1)k\right)}{2k} + \gamma q_1\left( (n+1)k\right) \notag \\
		&- \mathbf{N}(2\mathbf{u}_1(nk) -  \mathbf{u}_1\left( (n-1)k\right)) \cdot \mathbf{u}_1\left( (n+1)k\right) = \gamma + \varepsilon_{q_1, n},
	\end{align}
	
	\begin{equation}
		||\varepsilon_{u_j, n}||, |\varepsilon_{q_j, n}| \leq C(R_1, T) k^2, \quad \forall n, s.t. \, nk \in [0,T]. \label{eq:RHS_boundedness}
	\end{equation}
	\label{eq:discrete_continuous_trajectory}
\end{subequations}

Taking the difference between \eqref{eq:discrete_continuous_trajectory} and \eqref{eq:2_Novel_SAV_discretization} with initial conditions $\begin{bmatrix}
	\mathbf{v}^0 \\
	q^0
\end{bmatrix}$ and $\begin{bmatrix}
	\mathbf{v}^1 \\
	q^1
\end{bmatrix}$, 
denoting $\begin{bmatrix}
	\delta \mathbf{u}_j^n \\
	\delta q_j^n 
\end{bmatrix} = \begin{bmatrix}
	\mathbf{u}_j(nk) - \mathbf{u}_j^n \\
	q_j(nk) - q_j^n
\end{bmatrix}, j = 0, 1$, where $\begin{bmatrix}
	\mathbf{u}_j^n \\
	q_j^n 
\end{bmatrix}  = S^n_{k} \begin{bmatrix}
	u_j \\
	q_j
\end{bmatrix}$, $j = 0, 1$. 
Here, $S_{k}$ is the solution operator to the novel mr-SAV-BDF2 scheme  $\eqref{eq:2_Novel_SAV_discretization}$. We deduce
\begin{subequations}
	\begin{align}
		&\frac{3 \delta \mathbf{u}^{n+1}_j - 4 \delta \mathbf{u}_j^n + \delta \mathbf{u}_j^{n-1}}{2k} + A \delta \mathbf{u}^{n+1}_j + \delta q_j^{n+1} \cdot \mathbf{N}(2\mathbf{u}_j(nk) - \mathbf{u}_j((n-1)k)) \notag \\ 
		&+ q_j^{n+1} \left\{ \mathbf{N}(2\mathbf{u}_j(nk) - \mathbf{u}_j((n-1)k)) - \mathbf{N}(2\mathbf{u}_j^n - \mathbf{u}_j^{n-1}) \right\} = \varepsilon_{u_j, n}, \quad j = 0, 1 \label{eq:difference_u}
	\end{align}
	\begin{align}
		&\frac{3 \delta q_j^{n+1} - 4\delta q_j^{n} + \delta q_j^{n-1}}{2k} + \gamma \delta q^{n+1}_j - \mathbf{N}(2\mathbf{u}_j(nk) - \mathbf{u}_j((n-1)k)) \cdot \mathbf{u}_j((n+1)k) \notag \\
		&+ \mathbf{N}(2\mathbf{u}_j^n - \mathbf{u}_j^{n+1}) \cdot \mathbf{u}_j^{n+1} = \varepsilon_{q_j, n}, \quad j = 0, 1 \label{eq:difference_q}
	\end{align}
\end{subequations}

By uniform boundedness of discrete solutions (Proposition 2), continuous solution, and the local Lipschitz continuity of \textbf{N}, we deduce

\begin{subequations}
	\begin{equation}
		|\delta q_j^{n+1} \mathbf{N}(2\mathbf{u}_j(nk) - \mathbf{u}_j((n-1)k)| \leq C|\delta q^{n+1}_j|,
	\end{equation}
	\begin{align}
		|q_j^{n+1} \{ \mathbf{N}(2\mathbf{u}_j(nk) - \mathbf{u}_j((n-1)k)) &- \mathbf{N}(2\mathbf{u}_j^n - \mathbf{u}_j^{n-1}) \} \notag \\
		&\leq C \left( |\delta \mathbf{u}_j^n| + |\delta \mathbf{u}_j^{n-1}| \right),
	\end{align}
	\begin{align}
		|\mathbf{N}(2\mathbf{u}_j(nk) &- \mathbf{u}_j((n-1)k)) \cdot \mathbf{u}_j((n+1)k) - \mathbf{N}(2\mathbf{u}_j^n - \mathbf{u}_j^{n+1}) \cdot \mathbf{u}_j^{n+1}| \notag \\ 
		&\leq C \left( |\delta \mathbf{u}_j^{n+1}| + |\delta \mathbf{u}_j^{n}| + |\delta \mathbf{u}_j^{n-1}| \right), \quad j = 0,1.
	\end{align}
	\label{eq:Proposition_2_eq}
\end{subequations}

$\eqref{eq:difference_u} \cdot \delta \mathbf{u}^{n+1}_j + \eqref{eq:difference_q} \cdot \delta q_j^{n+1}$ and utilize \eqref{eq:Proposition_2_eq} and \eqref{eq:RHS_boundedness}. We deduce 

\begin{align*}
	&\frac{1}{2k} \{ ||\delta \mathbf{V}_j^{n+1}||^2_G - ||\delta \mathbf{V}_j^{n}||^2_G + ||\delta \mathbf{Q}_j^{n+1}||^2_G - ||\delta \mathbf{Q}_j^{n}||^2_G + \frac{1}{4} || \delta \mathbf{u}_j^{n+1} - 2\delta \mathbf{u}_j^{n} \\
	&+\delta \mathbf{u}_j^{n-1}||^2 + \frac{1}{4} | \delta q_j^{n+1} - 2q_j^{n} +\delta q_j^{n-1}|^2 \} + ||A^{\frac{1}{2}} \mathbf{u}^{n+1}_j ||^2 + \gamma |\delta q^{n+1}_j|^2 \\
	&\leq C|\delta q^{n+1}_j| \ ||\delta \mathbf{u}^{n+1}_j|| + C\left( ||\delta \mathbf{u}^{n}_j|| + ||\delta \mathbf{u}^{n-1}_j|| \right) ||\delta \mathbf{u}^{n+1}_j|| \\ 
	&+ C|\delta q^{n+1}_j| \left( ||\delta \mathbf{u}^{n+1}_j|| + ||\delta \mathbf{u}^{n}_j|| + ||\delta \mathbf{u}^{n-1}_j|| \right) + Ck^2 \left( ||\delta \mathbf{u}^{n}_j|| + |\delta q^{n}_j| \right).
\end{align*}

Applying the Cauchy-Schwartz, and utilizing the positivity of A and the equivalence of the G-norm and the standard norm on $\mathbb{H} = H \times H$ and $\mathbb{R}^2$, we deduce
\begin{align}
	&\frac{1}{2k} \left\{ ||\delta \mathbf{V}_j^{n+1}||^2_G - ||\delta \mathbf{V}_j^{n}||^2_G + ||\delta \mathbf{Q}_j^{n+1}||^2_G - ||\delta \mathbf{Q}_j^{n}||^2_G \right\} \notag \\
	&\leq C \left( ||\delta \mathbf{V}^{n+1}_j||_G^2 + ||\delta \mathbf{V}^{n}_j||_G^2 + ||\delta \mathbf{Q}^{n+1}_j||_G^2 + ||\delta \mathbf{Q}^{n}_j||_G^2 \right) + Ck^4.  
\end{align}
This implies, for $k < \frac{1}{2C}$,
\begin{align}
	||\delta \mathbf{V}_j^{n+1}||^2_G + ||\delta \mathbf{Q}_j^{n+1}||^2_G &\leq \left(\frac{1 + 2kC}{1-2kC}\right)^n \left( ||\delta \mathbf{V}^{1}_j||_G^2  + ||\delta \mathbf{Q}^{1}_j||_G^2 + \frac{Ck^4}{2} \right) \notag \\
	&\leq C(T, k) \left( ||\delta \mathbf{V}^{1}_j||_G^2  + ||\delta \mathbf{Q}^{1}_j||_G^2 + k^4 \right), \quad  j = 0, 1.
\end{align}

Thanks to Proposition 3, and $\eqref{Eq:proposition_3}$, we have 
$\begin{bmatrix}
	\delta\mathbf{V}_0^1 \\
	\delta \mathbf{Q}_0^1
\end{bmatrix} = \left[ S_q(k) \begin{bmatrix}
	\mathbf{v}^0\\
	q^0
\end{bmatrix} - \begin{bmatrix}
	\mathbf{v}^0 \\
	q^0
\end{bmatrix} \right] + \left[ \begin{bmatrix}
	\mathbf{v}^0 \\
	q^0 
\end{bmatrix} - \begin{bmatrix}
	\mathbf{v}^1 \\
	q^1
\end{bmatrix} \right]$. However, 

\begin{align*}
	\left|S_q(k) \begin{bmatrix}
		\mathbf{v}_0 \\
		q_0 
	\end{bmatrix} - \begin{bmatrix}
		\mathbf{v}_0 \\
		q_0
	\end{bmatrix} \right| &\leq Ck, \\
	\left| \begin{bmatrix}
		\mathbf{v}_1 \\
		q_1 
	\end{bmatrix} - \begin{bmatrix}
		\mathbf{v}_0 \\
		q_0
	\end{bmatrix} \right| &\leq Ck.
\end{align*}

Likewise 
$\delta \mathbf{V}_1^1 = \mathscr{O}(\delta t) = \mathscr{O}(k), \quad \delta \mathbf{Q}_1^1 = \mathscr{O}(\delta t) = \mathscr{O}(k)$. Therefore, we have to prove the following result.

\begin{proposition}
	$\forall T>0, \quad \exists C(T, R_1) > 0,$ s.t. 
	\begin{equation}
		||\delta \mathbf{V}_j^{n+1}||^2_G + ||\delta \mathbf{Q}_j^{n+1}||^2_G \leq C(T,R_1) k^2, \quad \forall (\mathbf{V}^0, \mathbf{Q}^0) \in \mathscr{A}_k, \quad nk \leq T.
	\end{equation}
\end{proposition}
\begin{remark} The result is not the expected second order for a second order scheme. This is due to the first order error committed at the first step.  If we choose the initial data so that $\delta \mathbf{V}_j^1 =  \mathscr{O}(k^2), \quad \delta \mathbf{Q}_j^1 =  \mathscr{O}(k^2), j=0,1$, the second order convergence follows.
\end{remark}

\section{Convergence of the attractor}

The purpose of this section is to show that the global attractors of the scheme, after taking appropriate projection, converge to that of the original model. Recall that the convergence of the global attractor is not automatic even for convergent numerical schemes \cite{SH1996}.

Recall that the classical result on convergence of attractor \cite{HLR1988} is not directly applicable since we do not have the convergence of the scheme in a neighborhood of the global attractor. See also \cite{Shen1989, Shen1990}. Likewise, the monoid approach proposed in \cite{HillSuli2000} is not directly applicable either since our scheme does not fall into the linear multi-step methods considered in that work.   Instead, we follow the approach of \cite{Wang2010} by focusing on the attractors of the discrete system and utilize the asymptotic consistency result as well as the convergence result that we derived in the previous two sections.

We first show that $\text{dist}\left( \mathscr{A}_k, \mathscr{A}_q \times \mathscr{A}_q \right) \underset{k \to 0}{\longrightarrow} 0$, where \\
$\mathscr{A}_q \times \mathscr{A}_q = \left\{ \left( \begin{bmatrix}
	\mathbf{u}\\
	q
\end{bmatrix}, \begin{bmatrix}
	\Tilde{\mathbf{u}} \\
	\Tilde{q}
\end{bmatrix} \right)  \Bigg|  \begin{bmatrix}
	\mathbf{u}\\
	q
\end{bmatrix}, \begin{bmatrix}
	\Tilde{\mathbf{u}} \\
	\Tilde{q}
\end{bmatrix} \in \mathscr{A}_q \right\},$\\ with $\mathscr{A}_q$ being the global attractor of the augmented system \eqref{eq:2_Novel_SAV_discretization}. 

Thanks to Proposition 3, $\mathscr{A}_k \subset B_{R_1} \subset \mathbb{H} \times \mathbb{R}^2, \forall k$.
By the fact that $\mathscr{A}_q \times \mathscr{A}_q$ is attracting for $\begin{bmatrix}
	S_q(t)\\
	S_q(t)
\end{bmatrix} = \mathbb{S}(t)$, we deduce, 
for any $\varepsilon >0, \exists T(\varepsilon, R_1) > 0, $ s.t. 
\begin{equation}
	\text{dist}\left(\mathbb{S}(t) \mathscr{A}_k, \mathscr{A}_q \times \mathscr{A}_q \right) \leq \frac{\varepsilon}{2}, \quad \forall t \geq T(\varepsilon, R_1).
\end{equation}

Let $n_k = \lfloor \frac{T(\varepsilon, R_1)}{k}\rfloor$. For any $ (\mathbf{V}^0, \mathbf{Q}^0) \in\mathscr{A}_k, \quad \mathbf{V}^0 = \begin{bmatrix}
	\mathbf{v}^0 \\
	\mathbf{v}^1
\end{bmatrix}, \quad \mathbf{Q}^0 = \begin{bmatrix}
	q^0 \\
	q^1
\end{bmatrix}, \quad \\ \exists (\Tilde{\mathbf{V}}^0, \Tilde{\mathbf{Q}}^0) \in \mathscr{A}_k, $ s.t. $(\mathbf{V}^0, \mathbf{Q}^0) = \mathbb{S}_k^{n_k} (\Tilde{\mathbf{V}}^0, \Tilde{\mathbf{Q}}^0) $ by the invariance of $\mathscr{A}_k$.
Hence
\begin{align*}
	\text{ dist}\left( (\mathbf{V}^0, \mathbf{Q}^0), \mathbb{S}(n_k k) \left[ \Tilde{\mathbf{V}}^0, \Tilde{\mathbf{Q}}^0 \right]\right) &= \text{dist}\left( \mathbb{S}_k^{n_k} \left[ \Tilde{\mathbf{V}}^0, \Tilde{\mathbf{Q}}^0 \right], \mathbb{S}(n_k k) \left[ \Tilde{\mathbf{V}}^0, \Tilde{\mathbf{Q}}^0 \right] \right), \\
	&\leq Ck \quad \text{(by Proposition 4).}
\end{align*}

Therefore, 
\begin{align*}
	\text{dist} \left( \left(\mathbf{V}^0, \mathbf{Q}^0 \right), \mathscr{A}_q \times \mathscr{A}_q  \right) &\leq \text{dist}\left( \left( \mathbf{V}^0, \mathbf{Q}^0 \right), \mathbb{S}(n_k k) \mathscr{A}_k \right) + \text{dist}\left( \mathbb{S}(n_k k) \mathscr{A}_k, \mathscr{A}_q \times \mathscr{A}_q  \right) \\
	&\le Ck+  \varepsilon.
\end{align*}
Hence, we deduce
$$\text{dist} \left( \mathscr{A}_k,\mathscr{A}_q \times \mathscr{A}_q  \right) \leq \varepsilon + Ck.$$
Letting $k \rightarrow 0$, we have 
$$\lim\limits_{k\rightarrow0}\text{dist}\left(\mathscr{A}_k, \mathscr{A}_q \times \mathscr{A}_q  \right) \leq \varepsilon.$$
Since $\varepsilon$ is arbitrary, we have derived

\begin{proposition}
	$\text{dist} \left(\mathscr{A}_k, \mathscr{A}_q \times \mathscr{A}_q  \right) \underset{k \to 0}{\longrightarrow} 0$. 
\end{proposition}
Notice that $\forall (\mathbf{V}^0, \mathbf{Q}^0) \in \mathscr{A}_k,$ 
$||\mathbf{v}^0 - \mathbf{v}^1|| \leq Ck, \quad |q^0 - q^1| \leq Ck$, by \eqref{Eq:proposition_3}.
Hence, the limit of $\mathscr{A}_k$ must live on the diagonal of the product space $\mathbb{H}\times \mathbb{R}^2$. Henceforth, we arrived at the following result.
\begin{theorem}
	$$ \text{dist}(\mathscr{A}_k, \mathds{1} \mathscr{A}_q) \underset{k \to 0}{\longrightarrow} 0, $$
	where $\mathds{1}\mathscr{A}_q = \left\{ \left( \begin{bmatrix}
		\mathbf{u}\\
		q
	\end{bmatrix}, \begin{bmatrix}
		\mathbf{u} \\
		q
	\end{bmatrix} \right) \in \mathscr{A}_q \times \mathscr{A}_q \right\}.$
	In particular, if we denote $\mathcal{P}_j$ as the projection onto the $j^{th}$ coordinate of the product space, $j=1,2$, and denoting $\mathcal{P}_u$ as the projection from $H\times \mathcal{R}^1$ to $H$, we have, when combined with proposition 2, 
	\begin{equation}
		\text{dist}(\mathcal{P}_u\mathcal{P}_j\mathscr{A}_k, \mathscr{A}) \underset{k \to 0}{\longrightarrow} 0,  j=1,2.
		\label{attractor_convergence}
	\end{equation}   
\end{theorem}

\section{Convergence of long-time statistics}

The convergence of the global attractor is an excellent indicator that the long time behavior of the scheme is closely related to that of the original underlying model, especially in terms of the geometry. However, even if two systems share the same global attractor could have completely different dynamics on the global attractor. 
For chaotic and/or turbulent systems, many of their physical properties are revealed statistically \cite{MW2006, Majda2016, LM1994,FMRT2001, F1995, SS2001}.
Therefore, it is of great importance to study whether the long time statistics of the scheme approximate that of the underlying model.
One of the important objects that characterizes the long time statistics is the invariant measures. A probability measure $\mu$ on the phase space is called invariant under the semigroup $S(t)$ if $\mu(S^{-1}(t)B)=\mu(B)$ for all measurable sets $B$. Our main goal in this section is to show that the invariant measures of the scheme, after taking appropriate marginal distribution, converge to some invariant measure of the underlying model.

Let $\mu_k \in \mathcal{P M}_k = \left\{ \text{invariant measure of } \mathbb{S}_k \text{ on } \mathbb{H} \times \mathbb{R}^2 \right\}$. 
Let $\mathcal{P}_j$ be the projection from $\mathbb{H} \times \mathbb{R}^2 \longrightarrow H \times R$ in the $j^{th}$ component. These projections induce  projections $\mathcal{P}_j^*$ on the space of probability measures. 
$$(\mathcal{P}_j^* \mu) (\varOmega) = \mu(P_j^{-1}(\varOmega)), \forall \varOmega \in \mathscr{B}({H} \times \mathbb{R}^1).$$ 

Let $ \Phi (\mathbf{u}, q)$ be a smooth test functional with compact support on $H \times \mathbb{R}^1$. Thanks to Proposition 2 and the fact that all invariant measures are supported on the global attractor \cite{Wang2009, FMRT2001}, we see that $\{\mu_k, k \geq 0\}$ form a weakly pre-compact set in the space of probability measures on $\mathbb{H} \times \mathbb{R}^2$. Hence, $\exists$ a subsequence, still denoted $\{\mu_k\}$ and a $\mathbf{\mu}_0 \in \mathcal{PM} (\mathbb{H} \times \mathbb{R}^2)$, s.t. 
$\mu_k \rightharpoonup \mu_0$, weakly in $\mathcal{PM}(\mathbb{H} \times \mathbb{R}^2)$.

Our goal is to show that $P_j^* \mu_k \rightharpoonup P_j^* \mu_0 \in \mathcal{PM}(S_q)$,  i.e., $P_j^* \mu_0$ is an invariant measure of system \eqref{eq:2_Novel_SAV_discretization}. 
We follow an argument similar to those presented in \cite{Wang2012}.
This is equivalent to showing, in the weak form of the invariance,
\begin{equation}
	\int_{H \times \mathbb{R}^1} \begin{bmatrix}
		- \mathbf{F} + A\mathbf{u} + q\mathbf{N}(\mathbf{u}) \\
		\gamma q - \gamma 
	\end{bmatrix} \cdot \Phi' \left(\begin{bmatrix}
		\mathbf{u}\\
		q
	\end{bmatrix} \right) d \left( P_j^* \mu_0 \right) \left(\begin{bmatrix}
		\mathbf{u}\\
		q
	\end{bmatrix} \right) = 0, \quad j = 1, 2,  \quad \forall \Phi. \label{eq:Convergnece_stats_goal}
\end{equation}

Without loss of generality, we assume $j = 1$. 

Let $ \mathbf{V} = \begin{bmatrix}
	\mathbf{u}_1\\
	\mathbf{u}_0
\end{bmatrix} \text{and } \mathbf{Q} = \begin{bmatrix}
	q_1\\
	q_0
\end{bmatrix}$.  Then the solution semigroup generated by the novel scheme \eqref{eq:2_Novel_SAV_discretization} on the product space gives us
$$\mathbb{S}_{k} \begin{bmatrix}
	\mathbf{V}\\
	\mathbf{Q}
\end{bmatrix}=  \begin{bmatrix}
	[\mathbf{u}_2,\mathbf{u}_1]^T\\
	[q_2, q_1]^T
\end{bmatrix}.
$$
Hence, the LHS of \eqref{eq:Convergnece_stats_goal} becomes

\begin{align*}
	LHS &= \int_{H \times \mathbb{R}^1} \begin{bmatrix}
		- \mathbf{F} + A\mathbf{u}_0 + q_0\mathbf{N}(\mathbf{u}_0) \\
		\gamma q_0 - \gamma 
	\end{bmatrix} \cdot \Phi' \left(\begin{bmatrix}
		\mathbf{u}_0\\
		q_0
	\end{bmatrix} \right) d P_1^* \mu_0 \\
	&=  \int_{\mathbb{H} \times \mathbb{R}^2}
	\begin{bmatrix}
		- \mathbf{F} + A\mathbf{u}_0 + q_0\mathbf{N}(\mathbf{u}_0) \\
		\gamma q_0 - \gamma 
	\end{bmatrix} \cdot \Phi' \left(\begin{bmatrix}
		\mathbf{u}_0\\
		q_0
	\end{bmatrix} \right) d  \mu_0 \\
	&= \lim_{k\rightarrow 0} \int_{\mathbb{H} \times \mathbb{R}^2} \begin{bmatrix}
		- \mathbf{F} + A\mathbf{u}_1 + q_1\mathbf{N}(\mathbf{u}_1) \\
		\gamma q_1 - \gamma 
	\end{bmatrix} \cdot \Phi' \left(\begin{bmatrix}
		\mathbf{u}_0\\
		q_0
	\end{bmatrix} \right) d \mu_k \left(\begin{bmatrix}
		\mathbf{V}\\
		\mathbf{Q}
	\end{bmatrix} \right) \\
	&= \lim_{k\rightarrow 0} \int_{\mathbb{H} \times \mathbb{R}^2} \begin{bmatrix}
		- \mathbf{F} + A\mathbf{u}_2 + q_2\mathbf{N}(2\mathbf{u}_1 - \mathbf{u}_0) \\
		\gamma q_2 - \gamma - \mathbf{N}(2\mathbf{u}_1 - \mathbf{u}_0) \cdot \mathbf{u}_2
	\end{bmatrix} \cdot \Phi' \left( \begin{bmatrix}
		\frac{3}{2} \mathbf{u}_1 - \frac{1}{2} \mathbf{u}_0 \\
		\frac{3}{2} q_1 - \frac{1}{2} q_0
	\end{bmatrix} \right) d \mu_k \left(\begin{bmatrix}
		\mathbf{V}\\
		\mathbf{Q}
	\end{bmatrix} \right)\\
	&\qquad  ( \text{Proposition 2\&3,  local Lipschitz and energy-conservation of \textbf{N}, and the smoothness and compact support of $\Phi'$}) \\
	&= \lim_{k\rightarrow 0} \int_{\mathbb{H} \times \mathbb{R}^2}  -\frac{1}{2k } \begin{bmatrix}
		3\mathbf{u}_2 - 4\mathbf{u}_1 + \mathbf{u}_0\\
		3q_2 - 4q_1 + q_0
	\end{bmatrix} \cdot \Phi' \left(\begin{bmatrix}
		\frac{3}{2} \mathbf{u}_1 - \frac{1}{2} \mathbf{u}_0 \\
		\frac{3}{2} q_1 - \frac{1}{2} q_0
	\end{bmatrix} \right) d \mu_k \left( \begin{bmatrix}
		\mathbf{V}\\
		\mathbf{Q}
	\end{bmatrix} \right) \\
	& \qquad ( \text{We have used scheme \eqref{eq:2_Novel_SAV_discretization}})\\
	&= - \lim_{k \rightarrow 0} \frac{1}{2k} \int_{\mathbb{H} \times \mathbb{R}^2} 
	\begin{bmatrix}
		\begin{bmatrix}
			3\\
			-1
		\end{bmatrix} \cdot \left( \begin{bmatrix}
			\mathbf{u}_2 \\
			\mathbf{u}_1
		\end{bmatrix}  - \begin{bmatrix}
			\mathbf{u}_1\\
			\mathbf{u}_0
		\end{bmatrix} \right) \\
		\begin{bmatrix}
			3\\
			-1
		\end{bmatrix} \cdot \left( \begin{bmatrix}
			q_2 \\
			q_1
		\end{bmatrix}  - \begin{bmatrix}
			q_1\\
			q_0
		\end{bmatrix} \right)
	\end{bmatrix} \cdot \Phi' \left( \begin{bmatrix}
		\frac{3}{2} \mathbf{u}_1 - \frac{1}{2} \mathbf{u}_0 \\
		\frac{3}{2} q_1 - \frac{1}{2} q_0
	\end{bmatrix} \right) \, d \mu_k \left(\begin{bmatrix}
		\mathbf{V}\\
		\mathbf{Q}
	\end{bmatrix} \right) \\
	&= \lim_{k \rightarrow 0} \frac{1}{2k} \int_{\mathbb{H} \times \mathbb{R}^2} 
	\begin{bmatrix}
		\begin{bmatrix}
			3\\
			-1
		\end{bmatrix} \cdot \left( \mathbb{S}_{k} \left(\begin{bmatrix}
			\mathbf{u}_1 \\
			\mathbf{u}_0
		\end{bmatrix}  , \begin{bmatrix}
			q_1\\
			q_0
		\end{bmatrix} \right)_u - \begin{bmatrix}
			\mathbf{u}_1 \\
			\mathbf{u}_0
		\end{bmatrix} \right)\\
		\begin{bmatrix}
			3\\
			-1
		\end{bmatrix} \cdot \left( \mathbb{S}_{k} \left(\begin{bmatrix}
			\mathbf{u}_1 \\
			\mathbf{u}_0
		\end{bmatrix}  , \begin{bmatrix}
			q_1 \\
			q_0
		\end{bmatrix} \right)_q - \begin{bmatrix}
			q_1 \\
			q_0
		\end{bmatrix}\right)
	\end{bmatrix} \cdot \Phi' \left(\begin{bmatrix}
		\frac{3}{2} \mathbf{u}_1 - \frac{1}{2} \mathbf{u}_0 \\
		\frac{3}{2} q_1 - \frac{1}{2} q_0
	\end{bmatrix} \right) d \mu_k \left( \begin{bmatrix}
		\mathbf{V} \\
		\mathbf{Q}
	\end{bmatrix} \right),\\
	&\qquad ( \left(\mathbb{S}_{k} \left[ \ \right] \right)_u \text{ denotes the u component}, \left(\mathbb{S}_{k} \left[ \ \right] \right)_q \text{ denotes the q component} ) \\
	&= \lim_{k\rightarrow0} \frac{1}{k} \int_{\mathbb{H} \times \mathbb{R}^2} 
	\left(
	\mathbb{S}_{k} 
	\begin{bmatrix}
		\mathbf{V}\\
		\mathbf{Q}
	\end{bmatrix} \cdot 
	\begin{bmatrix}
		-\frac{3}{2} \\
		-\frac{1}{2}
	\end{bmatrix}
	-  \begin{bmatrix}
		\mathbf{V} \\
		\mathbf{Q}
	\end{bmatrix} \cdot 
	\begin{bmatrix}
		\frac{3}{2} \\
		-\frac{1}{2}
	\end{bmatrix} 
	\right)
	\cdot
	\Phi' \left(
	\begin{bmatrix}
		\mathbf{V} \\
		\mathbf{Q}
	\end{bmatrix} \cdot \begin{bmatrix}
		\frac{3}{2} \\
		-\frac{1}{2}
	\end{bmatrix}
	\right)
	d \mu_k \left(\begin{bmatrix}
		\mathbf{V} \\
		\mathbf{Q}
	\end{bmatrix} \right) \\
	&= \lim_{k\rightarrow0} \frac{1}{k} \int_{\mathbb{H} \times \mathbb{R}^2}
	\left( \Phi \left(\mathbb{S}_{k} \begin{bmatrix}
		\mathbf{V}\\
		\mathbf{Q}
	\end{bmatrix} \cdot \begin{bmatrix}
		-\frac{3}{2} \\
		-\frac{1}{2}
	\end{bmatrix}  \right)-  \Phi \left(  \begin{bmatrix}
		\mathbf{V} \\
		\mathbf{Q}
	\end{bmatrix} \cdot \begin{bmatrix}
		\frac{3}{2} \\
		-\frac{1}{2}
	\end{bmatrix} \right) \right)
	d \mu_k \left( \begin{bmatrix}
		\mathbf{V} \\
		\mathbf{Q}
	\end{bmatrix} \right)\\
	&\qquad ( \text{ Taylor expansion and Proposition 3} )\\
	&= 0 , \\
	&\qquad( \text{Since } \mu_k \text{ is invariant under }\mathbb{S}_k ).
\end{align*}
Hence, we have proved the following result.
\begin{theorem}
	Let $\mu_k$ be an invariant measure of the semigroup on the product space $\mathbb{H}\times\mathbb{R}^2$ generated by the novel mean-reverting-SAV-BDF2 scheme \eqref{eq:2_Novel_SAV_discretization} with time step $k$. Let $\mathcal{P}_j$ be the projection from $\mathbb{H}\times\mathbb{R}^2$ onto its $j^{th}$ component, and let $\mathcal{P}^*_j$ be the induced marginal distribution in the $j^{th}$ component. Then $\{\mu_j\}$ contains a subsequence, still denoted $\{\mu_j\}$,  that weakly converges to a probability measure $\mu_0$ on $\mathbb{H}\times\mathbb{R}^2$. Moreover, $\mathcal{P}^*_j\mu_0$ is an invariant measure of $S_q$, the solution semi-group of the extended system. In addition, if we denote $\mathcal{P}_u$ be the projection from $H\times \mathbb{R}^1$ to the $\mathbf{u} $ component, and let $\mathcal{P}^*_u\mu$ be the marginal distribution of $\mu$ in the $\mathbf{u}$ component, then 
	\begin{equation}
		\mathcal{P}^*_u\mathcal{P}^*_j\mu_k\rightharpoonup  \mathcal{P}^*_u\mathcal{P}^*_j\mu_0\in \mathcal{IM}(S),
		\label{IM_convergence}
	\end{equation}
	where $\mathcal{IM}(S)$ denotes the set of invariant measures of the solution semigroup $S$ to the original model.
	
\end{theorem}

\section{Application to Lorenz 96}

In this section, we employ the novel mr-SAV-BDF2 scheme to study the long time statistical properties of the following damped and driven Lorenz 96 model 
\cite{lorenz1996, lorenz1998,KloedenLorenz1990,MW2006,MAG2005}. 
{
	This is a toy model for atmospheric motion  that enjoys some of the general properties of geophysical models, namely energy-preserving advection, damping and forcing \cite{MAG2005,lorenz1996, lorenz1998}.}
\begin{equation}
	\frac{du_j}{dt} = \left(u_{j+1} - u_{j-2} \right) u_{j-1} - u_j + F, 
	\label{eq:L96_model_j}
\end{equation}
where $j$ is the index for spatial locations, and $F$ is the forcing term. We assume periodic boundary condition, i.e., $u_{-1} = u_{J-1}, u_0 = u_{J},  u_{1} = u_{J+1}$. 
For simplicity we will focus on the case of $J=5$ and $F=-12$ as this case already exhibits strong chaotic behavior.
We will focus on the long time distribution using the long-time simulation of the system, i.e., the so-called indirect approach \cite{SS2001}.

It is known that the limit of the long time averages of either the original model or the numerical scheme correspond to invariant measures of the model or the scheme \cite{Wang2009, Wang2010, FMRT2001}. Therefore, the distribution of the trajectory of each coordinate in $H$ of the system converges to the corresponding marginal distribution of an invariant measure of the  the original model or the novel numerical scheme. Hence, the normalized long time histogram of each coordinate of the system should converge to an equilibrium as the time interval approaches infinity. However, the convergence results from the literature or the previous section do not provide the rate of convergence in terms of the length of the time interval, or the size of the time-step. Hence we resort to numerics. 

Our numerical experiments verify the second order accuracy of the novel BDF2-FSAV scheme on any finite time interval, as expected. However, the prefactor is not small. For $\delta t=k=2^{-14}$, the relative error is of the order of $10^{-6}$ on the time interval $[0,1]$  using a numerical truth generated by the scheme with a very small time step ($k=2^{-23}$) and random initial data drawn from $[-15,15]^5$.  Nevertheless, the error quickly grows to the order of $10^{-3}$ on the time interval $[0,5]$, consistent with the expected chaotic behavior of the system. This suggests $10^{-3}$ as a threshold for $\delta t=2^{-14}$.
For the numerical experiments we have used the parameters, $F_j = -12 \ \forall j \in J, \gamma = 1000$. 
We report the result for the first location ($j=1$) as the long time statistics are expected to be similar due to the rotation invariance of the model. 
We have used the Jensen-Shannon (JS) entropy/distance/divergence/information radius, and the total variation distance to investigate the statistics. 
Recall the JS entropy between two discrete probability distributions $P$ and $Q$ on the same space are defined as 
\begin{equation}
	D_{JS}(P||Q) = \frac{1}{2} \left( D_{KL}(P||\frac{P+Q}{2}) + D_{KL}(Q||\frac{P+Q}{2}) \right)
\end{equation}
where
$$D_{KL}(P||Q) = \sum_{i} P(i) \log \left(\frac{P(i)}{Q(i)} \right).$$

\subsection{The convergence in terms of the terminal time}
Although the long time averages over time intervals $[0, T]$ are expected to converge as $T\rightarrow\infty$ \cite{FMRT2001, Wang2009, Wang2010}, no rate of convergence is known.
Therefore, we resort to numerics.
To investigate the impact of the length of the time interval we have used various time lengths against our long-time ``truth'' over the interval $[0, 2\times 10^6]$. Table  \ref{tab:Terminal_Time_test_j=1} clearly indicates a half order convergence rate in the total variation distance of the difference between the numerical  truth and the long time statistics over the time interval $[0, T]$.  If this half-order convergence result holds for all times, this would imply that there is no exponential mixing for this model in this strongly chaotic parameter regime, unlike the stochastic models \cite{E2001, SS2001}. In addition, this would also imply that one needs to perform simulation on a time interval of the order of $10^6$ time units in order to get the first three digits of the distribution correct.
This highlights the challenge in studying the climate of deterministic systems using indirect approach \cite{SS2001}.
Table \ref{tab:Terminal_Time_test_j=1} also shows first order convergence of the long-time statistics using the JS entropy/distance. The first order convergence of the JS distance/entropy is consistent with the half-order convergence of the total variation distance via the so-called Pinsker's inequality.\footnote{$\|P-Q\|_{TV}\le \sqrt{\frac12 D_{KL}(P\|Q)}$} 
%
\begin{table}[h]
	\begin{center}
		\begin{tabular}{||c c c c c||} 
			\hline
			$T$ & JS & Order &  TV  & Order \\[0.5ex] 
			\hline\hline
			100 & 3.3779e-02 &  & 3.5182e-01 &  \\ 
			\hline
			200 & 1.7769e-02 & 0.9268 & 2.5147e-01 & 0.4845 \\ 
			\hline
			400 & 9.2915e-03 & 0.9354 & 1.8068e-01 & 0.4769 \\ 
			\hline
			800 & 4.9389e-03 & 0.9117 & 1.3457e-01 & 0.4251 \\ 
			\hline
			1600 & 2.5206e-03 & 0.9704 & 9.4286e-02 & 0.5132 \\ 
			\hline
			3200 & 1.3782e-03 & 0.8710 & 6.8695e-02 & 0.4568 \\ 
			\hline
			6400 & 6.8670e-04 & 1.0050 & 4.7947e-02 & 0.5188 \\ 
			\hline
			12800 & 3.5437e-04 & 0.9544 & 3.4542e-02 & 0.4731 \\ 
			\hline
			25600 & 1.8230e-04 & 0.9589 & 2.4575e-02 & 0.4912 \\ 
			\hline
			51200 & 9.4037e-05 & 0.9550 & 1.7788e-02 & 0.4663 \\ 
			\hline
			102400 & 4.8429e-05 & 0.9574 & 1.2560e-02 & 0.5021 \\ 
			\hline
			204800 & 2.5748e-05 & 0.9114 & 9.0831e-03 & 0.4676 \\ 
			\hline
			409600 & 1.4058e-05 & 0.8731 & 6.6435e-03 & 0.4512 \\ 
			\hline
		\end{tabular}
		\caption{The JS divergence score and TV distance for different terminal times, $T$, against the reference solution over the interval $[0, 2\times 10^6]$, using $\delta t = 2^{-14}$ and 512,000 bins for the $j=1$ spatial location.}
		\label{tab:Terminal_Time_test_j=1}
	\end{center}
\end{table}


\subsection{The impact of the number of bins}
Although the distribution of the long-time trajectory is a distribution on a continuous space, we will have to use a finite probability to approximate the continuous space for our numerics. In the case of a scalar variable, such as the first component of the system, we partition the range of the variable into finitely many equal sized bins. 
To investigate the impact of the number of bins utilized to approximate the distribution of the statistical equilibrium, we compare the distribution of the signal over $[0, 2\times 10^6]$.
We use the moving average of the distribution on the interval $[0,2000000]$ with 512,000 bins as the reference solution for table \ref{tab:Bins_test_complete_statistics_Moving_Av}. 
\begin{table}[h]
	\begin{center}
		\begin{tabular}{||c c c||} 
			\hline
			N  & JS & TV  \\ [0.5ex] 
			\hline\hline
			125     & 7.3332e-05 & 1.5215e-02  \\
			\hline
			250     & 2.3680e-05 & 8.6328e-03  \\
			\hline
			500     & 1.0443e-05 & 5.6382e-03  \\
			\hline
			1000    & 6.5550e-06 & 4.3911e-03  \\
			\hline
			2000   & 4.6104e-06 & 3.6397e-03  \\
			\hline
			4000    & 3.1968e-06 & 3.0366e-03  \\
			\hline
			8000    & 1.9741e-06 & 2.3885e-03  \\
			\hline
			16000   & 9.9875e-07 & 1.6700e-03  \\
			\hline
			32000   & 4.0511e-07 & 1.0815e-03  \\
			\hline
			64000   & 1.6949e-07 & 7.1593e-04  \\
			\hline
			128000  & 7.2914e-08 & 4.7989e-04  \\
			\hline
			256000  & 2.7627e-08 & 3.0521e-04  \\
			\hline
		\end{tabular}
		\caption{The JS scores and the TV distance for different number of bins, $N$, against the reference solution with 512,000 bins for the $j=1$ spatial location on the interval $[0,2\times 10^6]$.}
		\label{tab:Bins_test_complete_statistics_Moving_Av}
	\end{center}
\end{table}
We observe that 64K bin leads to an error less than $10^{-3}$ the threshold value. Therefore, we use 64K or more  bins in subsequent experiments.


%
\subsection{The impact of the time-step size}
To investigate the impact of the time-step size, we utilize a numerical truth generated with $k=2^{-17}$ over $[0,1\times 10^5]$ as the reference and vary the time step.  The following table \ref{tab:Varying_dt_test_j=1_T=100000} indicates that the error saturates  at $\delta=2^{-14}$. Note that with half-order convergence rate and terminal time of $10^5$, the distance from this long-time statistics to the true equilibrium is expected to be of the order of $10^{-2}$, in agreement with the saturation value reported in this table. This also partially justifies our usage of $\delta t=2^{-14}$ for even longer simulations. 
%
\begin{table}[h]
	\begin{center}
		\begin{tabular}{||c c c c c||} 
			\hline
			$\delta t$ & JS & Order & TV & Order \\[0.5ex] 
			\hline\hline
			$2^{-9}$ & 1.0353e-03 &  & 6.0927e-02 &  \\ 
			\hline
			$2^{-10}$ & 5.3745e-04 & 0.9458 & 4.3598e-02 & 0.4828 \\ 
			\hline
			$2^{-11}$ & 2.8628e-04 & 0.9087 & 3.1427e-02 & 0.4723 \\ 
			\hline
			$2^{-12}$ & 1.6479e-04 & 0.7968 & 2.3609e-02 & 0.4127 \\ 
			\hline
			$2^{-13}$ & 1.0995e-04 & 0.5838 & 1.9006e-02 & 0.3129 \\ 
			\hline
			$2^{-14}$ & 6.9974e-05 & 0.6520 & 1.4977e-02 & 0.3437 \\ 
			\hline
			$2^{-15}$ & 5.8395e-05 & 0.2610 & 1.3224e-02 & 0.1796 \\ 
			\hline
			$2^{-16}$ & 5.1675e-05 & 0.1764 & 1.2246e-02 & 0.1108 \\ 
			\hline
		\end{tabular}
		\caption{The JS divergence score and TV distance for different time steps, $\delta t$, over the time interval, $[0,100000]$ using 512,000 bins against the reference solution using $\delta t = 2^{-17}$ for the $j=1$ spatial location.}
		\label{tab:Varying_dt_test_j=1_T=100000}
	\end{center}
\end{table}

\subsection{The impact of the initial data}
To investigate the effect of the initial data we use two random initial conditions so that each component of the initial data is taken from a uniform distribution over $ [-15, 15]^5$, and the second one being $ 5\%$ random perturbation of the first in each component. We have varied the time intervals $[0, T]$, with $T$ being the end time of the interval. We use 64,000 bins for the investigation.
\begin{table}[h]
	\begin{center}
		\begin{tabular}{||c c c c c ||} 
			\hline
			T &  JS & Order &TV  & Order \\ [0.5ex] 
			\hline\hline
			100   & 2.6358e-02 &  & 3.2327e-01 &\\ 
			\hline
			200   &  1.2607e-02 & 1.0640 & 2.1569e-01 & 0.5838\\ 
			\hline
			400   & 8.7577e-03 & 0.5256 & 1.6356e-01 & 0.3991  \\ 
			\hline
			800   & 3.8982e-03 & 1.1677  & 1.1394e-01 & 0.5215 \\ 
			\hline
			1600  &  2.1966e-03 & 0.8275 & 8.1951e-02 & 0.4754 \\ 
			\hline
			3200  &  1.0673e-03 & 1.0413 & 5.8923e-02 & 0.4759 \\ 
			\hline
			6400  &  6.0011e-04 & 0.8307 & 4.3710e-02 & 0.4309\\ 
			\hline
			12800 & 2.7423e-04 & 1.1298  & 2.9404e-02 & 0.5720\\ 
			\hline
			25600 & 1.2791e-04 & 1.1003 & 1.9865e-02 & 0.5658\\ 
			\hline
		\end{tabular}
		\caption{The JS and TV distances and their orders for different terminal times, $T$, using the random and perturbed initial conditions with $64,000$ bins for the $j=1$ spatial location.}
		\label{tab:Initial_data_64k_bins_Norms_1_j=1}
	\end{center}
\end{table}
We observe first order convergence in the JS distance and half order convergence in TV similar to the convergence rate in the long-time statistics.
The table indicates two randomly chosen initial data lead to similar similar statistics over the long-time interval $[0, 25600]$ with the difference of the order of $10^{-2}$, consistent with the error from the true statistical equilibrium based on half-order convergence rate.
Tests with other randomly generated initial data yield qualitatively very similar results in our numerical experiments.

\subsection{Comparison to the first order scheme}

It is easy to have a highly efficient first order scheme that is able to capture the long time statistics via a combination of backward Euler, forced and damped SAV, and IMEX in exactly the same way as the second order one. 
\begin{align*}
	B^n &= 1 + \frac{\delta t^2}{1+ \gamma\delta t} \left( \mathbf{N}(\mathbf{u}^n) \cdot \left[ I + \delta t A \right]^{-1} \mathbf{N}(\mathbf{u}^n) \right) \\
	q^{n+1} &= \frac{1}{B^n \left(1+ \gamma\delta t \right)} \left( \delta t \gamma + q^n + \delta t \mathbf{N}(\mathbf{u}^{n}) \cdot \left(  \left[ I + \delta t A \right]^{-1} \left( \delta t \mathbf{F}^{n+1} + \mathbf{u}^n \right) \right) \right) \\
	\mathbf{u}^{n+1} &= \left[ I + \delta t A \right]^{-1} \left( \delta t \mathbf{F}^{n+1} + \mathbf{u}^n - \delta t q^{n+1} \mathbf{N}(\mathbf{u}^n)  \right)
\end{align*}
We naturally wonder whether the second order scheme performs better than the first order scheme. While the  rate of convergence in terms of the terminal time are the same, i.e.,  first order in the Jensen-Shannon distance, and half-order in TV, the second order scheme requires less time to achieve two digits of accuracy in both the mean and the variance as evident from the following table.
We use the reference solution on the time interval $[0,  4\times 10^{6}]$ generated by the second order scheme with $\delta t=2^{-14}$ as the numerical truth.
For the first component, the mean is $-2.30785305840738$, and the variance is $22.3539129577942$.
For the first order scheme , we use a smaller time step of $\delta t = 2^{-18}$.
\begin{table}[h]
	\begin{center}
		\begin{tabular}{||c | c c | c c||} 
			\hline
			Relative error & Mean $1^{st}$ order & Mean $2^{nd}$ order & Var $1^{st}$ order & Var $2^{nd}$ order \\ [0.5ex] 
			\hline\hline
			1\% & 42400 & 11700 & 12500 & 4000 \\ 
			\hline
		\end{tabular}
		\caption{The Terminal Time $T$ needed for the mean and variance to stay below the threshold of $1\%$ for the relative error. }
		\label{tab:Mean_Variance_Comparison_FO_SO}
	\end{center}
\end{table}
The result in this table illustrates the advantage of the second order scheme.

\section{Concluding remarks}
We have proposed a highly efficient second-order mean-reverting-SAV-BDF2-based (mr-SAV-BDF2) numerical scheme for a class of finite dimensional nonlinear models. The scheme is unconditional stable, enjoys a uniform-in-time  bound for arbitrary initial data, and is able to capture the long-time statistics of the underlying model under appropriate assumptions. These assumptions are satisfied by a large family of geophysical models,  including the Lorenz 96 (L96) system. 

Our numerics on the 5-mode L96 with a moderate forcing $F=-12$ suggests that it takes a very long time for the system to reach statistical equilibrium. 
Our numerics imply a half-order convergence rate,  highlighting the need of extremely long simulations in order to control the total variation norm of the error within $1\%$.
This observation underscores the challenges associated with studying climate dynamics and climate change of geophysical models.

The application of the novel  mean-reverting-SAV-BDF2 scheme to infinite dimensional models, such as the two-dimensional Navier-Stokes equations,  presents additional challenges, 
{
	as the required compactness in the Lax-type criteria no longer follows from boundedness, among others factors.} Nevertheless, the uniform-in-time bounds can still be established for such systems, see \cite{HW2024} for the case of 2D NSE.  Challenges such as the convergence of invariant measures will be addressed in  future works.

	\section*{Acknowledgement}
The work of Han is supported by the National Science Foundation grant DMS-2310340. The work of Wang is supported by the National Natural Science Foundation of China grant 12271237 as well as the Gary Havener Endowment. The authors are listed in alphabetic order, not according to the contributions.

\end{document}